\documentclass{daj}

\usepackage{AKstyle}

\usepackage{comment}

\newcommand\faS{Y}
\newcommand\xiS{X}

\dajAUTHORdetails{%
  title = {Beyond Expansion IV: Traces of Thin Semigroups}, 
  author = {Jean Bourgain and Alex Kontorovich},
  keywords = {thin groups, affine sieve, additive energy},
}   

\dajEDITORdetails{%
   year={2018},
   number={6},
   received={25 April 2017},   
   revised={30 March 2018},    
   published={10 April 2018},  
   doi={10.19086/da.3471},       
}   

\begin{document}

\begin{frontmatter}[classification=text]


\author[jb]{Jean Bourgain
\thanks{Supported by NSF grant DMS-1301619.}}

\author[ak]{Alex Kontorovich
\thanks{Supported by
an NSF CAREER grant DMS-1455705, an
NSF FRG grant DMS-1463940, an Alfred P. Sloan Research Fellowship, and a
BSF grant. Some of this work was carried out thanks to support from a Yale Junior Faculty Fellowship and the Ellentuck Fund at IAS.}}

\begin{abstract}
This paper constitutes Part IV in our study of particular instances of the Affine Sieve, producing 
level%
s of distribution beyond
those
 attainable from expansion alone. Motivated by McMullen's 
 Arithmetic Chaos
 Conjecture 
regarding
low-lying closed geodesics on the modular surface defined over a given number field, we study the set of traces for certain sub-semi-groups of $\SL_{2}(\Z)$ corresponding to absolutely Diophantine numbers (see \secref{sec:arithC}). In particular, we are concerned with the level of distribution for this set. While the standard Affine Sieve procedure, combined with Bourgain-Gamburd-Sarnak's resonance-free region for the resolvent of a ``congruence'' transfer operator, produces {\it some} exponent of distribution $\ga>0$, we are able to produce the exponent $\ga=1/3-\vep$.
This recovers unconditionally the same exponent as what one would obtain under a Ramanujan-type conjecture for thin groups.
A key ingredient, of independent interest, is a bound on the additive energy of $\SL_2(\mathbb Z)$.
\end{abstract}
\end{frontmatter}


\section{Introduction}\label{sec:intro}

In this paper, we reformulate McMullen's (Classical) Arithmetic Chaos Conjecture (see \conjref{eq:ArithC}) as a local-global problem for the set of traces in certain thin semigroups, see \conjref{conj:locGlob}. 
Our main goal is to make some partial progress towards this conjecture by establishing strong levels of distribution for this trace set, see \secref{sec:rough}.

\subsection{Low-Lying Closed Geodesics 
With Fixed Discriminant}\label{sec:low} 

This paper is
motivated by the
 study
 of 
long closed geodesics on the modular surface defined over a given number field, which do not have high excursions into the cusp.
Let us make this precise.

To set notation, 
let $\bH$ denote the 
upper half plane, 
and let
$$
\cX=T^{1}(\SL_{2}(\Z)\bk\bH)\cong\SL_{2}(\Z)\bk\SL_{2}(\R)
$$ 
be the unit tangent
bundle of the modular surface.
A closed geodesic $\g$ on $\cX$ corresponds to 
 a hyperbolic matrix
 $M\in\SL_{2}(\Z)$ (more precisely its conjugacy class).
Let $\ga_{M}\in\dd\bH
$ be 
one
of  the two 
fixed points of $M$, the other 
being its Galois conjugate $\overline{\ga_{M}}$; then $\g
$ is the projection mod $\SL_{2}(\Z)$ of the geodesic connecting $\ga_{M}$ and $\overline{\ga_{M}}$.
We will say that $\g
$ is defined over the (real quadratic) field $K
=\Q(\ga_{M})$.
Let 
%
 $\gD_{M}$ 
 be
 the discriminant of $K
 $; 
 this 
 number
 is (up to factors of $4$) the square-free part of 
\be\label{eq:DMis}
D_{M}:=(\tr M)^{2}-4.
\ee

To study excursions into the cusp, let $\sY({\g})
$ denote 
  the largest imaginary part of $\g
  $ in the standard upper-half plane fundamental domain for 
the modular surface. 
Given a ``height'' $C>1$, we say that the closed geodesic $\g$ is {\bf low-lying} (of height $C$) if $\sY({\g})<C$. 
By the well-known connection  \cite{Humbert1916, Artin1924, Series1985} between continued fractions and the cutting sequence of the geodesic flow on $\cX$, 
the  condition that $\g$ be  low-lying can be 
reformulated as a
 Diophantine property on the fixed point $\ga_{M}$ of $M$, as follows.
Write the (eventually periodic) continued fraction expansion  
$$
\ga_{M}=
a_{0}+\cfrac{1}{a_{1}+\cfrac{1}{a_{2}+\ddots}}
=
[a_{0},a_{1},a_{2},\dots]
,
$$
as usual,
where
the numbers $a_{j}$ are called {partial quotients}. Given any 
$A\ge1$, we say that $\ga_{M}
$ is {\bf Diophantine} (of height $A$)  if all its partial quotients $a_{j}$ are 
bounded by
$A$. Then $\ga_{M}$ being Diophantine of height $A$ is essentially equivalent to $\g$ being low-lying of height $C=C(A)$.  

\begin{question}\label{ques:1}
{\rm Given a 
real quadratic  field $K$
and a height $C$,
can one 
find 
longer and longer
primitive
 closed geodesics defined over
$K$
 which are low-lying of height $C$? 
Equivalently, given a fixed fundamental discriminant $\gD>0$ and a height $A\ge1$, we wish to find larger and larger (non-conjugate, primitive, hyperbolic) matrices 
$M$ so that their fixed points $\ga_{M}$ 
are
Diophantine of height $A$, 
and so that
$t=\tr(M)$ solves the Pell equation $t^{2}-\gD s^{2}=4$; cf. \eqref{eq:DMis}.
If solutions exist, how rare/ubiquitous are they?
}
\end{question}

\subsection{Arithmetic Chaos}\label{sec:arithC}

On one hand, the answer to Question \ref{ques:1} 
is,
on average, negative. Indeed, 
generic long closed geodesics
equidistribute, so must
 have 
 high excursions into the cusp.
On the other hand, McMullen's (Classical) Arithmetic Chaos Conjecture (see \cite{McMullen2009, McMullenNotes} for the dynamical perspective and origin of this problem) predicts that solutions exist, and 
moreover have
positive entropy:
\begin{conj}[Arithmetic Chaos \cite{McMullenNotes}]\label{conj:McM}
There exists an absolute height $A\ge2$ so that, for any fixed real quadratic field $K$, 
the cardinality of the set
\be\label{eq:ArithC}
\bigg\{
%
[\overline{a_{0},a_{1},\dots,a_{\ell}}]\in K 
\quad :\quad
1\le a_{j}\le A
\bigg\}
\ee
grows exponentially, as $
\ell
\to\infty$. 
\end{conj}

\begin{rmk}
Though we have stated the conjecture with some absolute height $A$,
McMullen formulated this problem with $A=2$ (of course $A=1$ only produces the golden mean).
He further suggested
 it should also hold whenever the corresponding growth exponent exceeds $1/2$, see \rmkref{rmk:genAC}. 
\end{rmk}

\begin{rmk}
As pointed out to us by McMullen, one can also formulate a $\GL_{n}(\Z)$ version of Arithmetic Chaos by strengthening \cite[Conjecture 1.7 (3)]{McMullen2009} so as to postulate exponential growth
of periodic points 
instead of just infinitude.
\end{rmk}

It is not currently known whether the following much weaker statement is true: 
for some $A$ and every $K$,
the cardinality of the set \eqref{eq:ArithC} is unbounded. Even worse, it is not known whether \eqref{eq:ArithC} is eventually  non-empty,
that is, 
whether 
there exists an $A\ge2$ so that any $K$ contains at least one element which is Diophantine of height $A$.
%

Some progress towards  \conjref{conj:McM} appears in 
 \cite{Woods1978, Wilson1980, McMullen2009, Mercat2012},
where  special periodic patterns of partial quotients are constructed to lie in certain prescribed real quadratic fields.
These results prove that for any $K$, there exists an $A=A(K)$ so that the cardinality of \eqref{eq:ArithC} is unbounded with $
\ell
$; 
but exponential growth is not known in a single case.

In light of this conjecture, we call a number {\bf absolutely Diophantine} if it is Diophantine of height $A$ for some absolute constant $A\ge1$. That is, when we speak of a number being absolutely Diophantine, the height $A$ is 
fixed in advance.
In the next subsection, we describe a certain ``local-global'' conjecture which has the Arithmetic Chaos Conjecture as a 
consequence.
 

\subsection{The Local-Global Conjecture}

First we need some more notation.
Consider  a 
finite subset $\cA\subset\N$, which we call an {\bf alphabet}, and let 
$$
\fC_{\cA}=\{[a_{0},a_{1},\dots]:a_{j}\in\cA\}
$$ 
denote the 
set of all $\ga\in\R$ 
with
all partial quotients in $
\cA$. 
If $\max\cA\le A$ for an absolute constant $A$, then
every $\ga\in\fC_{\cA}$ is clearly absolutely Diophantine. 
Assuming $2\le|\cA|<\infty$, each $\fC_{\cA}$ is a Cantor set of some Hausdorff dimension 
$$
0<\gd_{\cA}<1
,
$$
and 
by choosing $\cA$ appropriately
(for example, $\cA=\{1,2,\dots, A\}$ with $A$ large), one can make $\gd_{\cA}$ 
arbitrarily close to $1$
\cite{Hensley1992}.

It is easy to see that the matrix
\be\label{eq:Majs}
M=
\mattwo{a_{0}}110
\mattwo{a_{1}}110
\cdots
\mattwo{a_{\ell}}110
\ee
(of determinant $\pm1$)
fixes the quadratic irrational
$$
\ga_{M}=[\overline{a_{0},a_{1},\dots,a_{\ell}}]
,
$$
so we introduce the semi-group\footnote{The superscript $+$ in \eqref{eq:cGcAdef} denotes generation as a semigroup, that is, no inverses.}
\be\label{eq:cGcAdef}
\cG_{\cA}:=
\<
\mattwo a110
\ 
:
\ 
a\in\cA
\>^{+}
\quad\subset\quad\GL_{2}(\Z)
\ee
of all such matrices whose fixed points $\ga_{M}$ lie in $\fC_{\cA}$.
Preferring to work in $\SL_{2}$, we immediately pass to the even-length (determinant-one) sub-semi-group
\be\label{eq:GcAis}
\G_{\cA}
:=
\cG_{\cA}\cap\SL_{2}(\Z)
,
\ee
which is (finitely) generated by the products $\mattwos a110\cdot\mattwos b110$, for $a,b\in\cA$.

Having 
accounted for
the ``low-lying'' (or Diophantine) criterion, we must study the discriminants, or what is essentially the same, the set 
$$
\cT_{\cA}:=\{\tr M\ 
: \ 
M\in\G_{\cA}\}\quad\subset\quad\Z
$$
of traces in $\G_{\cA}$.
Borrowing language from Hilbert's 11$^{th}$ problem on numbers represented by quadratic forms, we call an integer $t$ {\bf admissible} (for the alphabet $\cA$) if for every $q\ge1$, 
$$
t\in\cT_{\cA}(\mod q)
,
$$
that is, if $t$ passes all finite local obstructions.
\begin{rmk}\label{rmk:StrongApp}
If $\{1,2\}\subseteq\cA$, then, allowing inverses, the {\it group} $\<\G_{\cA}\>$ generated by the semigroup $\G_{\cA}$ is 
all of 
$
\SL_{2}(\Z),
$ 
and hence every integer is admissible.
In general, 
 Strong Approximation \cite{MatthewsVasersteinWeisfeiler1984} shows that admissibility can be checked using a single modulus $q(\cA)$.
\end{rmk}
 We say $t$ is {\bf represented} if $t\in\cT_{\cA}$,
and let
$\cM_{\cA}(t)$
 denote 
 its {\bf multiplicity}, 
$$
\cM_{\cA}(t):=\#\{M\in\G_{\cA}:\tr M=t\}.
$$
Since the entries of $\G_{\cA}$ are all positive, the multiplicity is always finite. 

The following conjecture
seems plausible.

\begin{conj}[Local-Global Conjecture for 
Traces]\label{conj:locGlob}
Let $\cA$ be an alphabet for which
the dimension $\gd_{\cA}$ exceeds $1/2$. 
Then the set $\cT_{\cA}$ of traces contains every sufficiently large admissible integer. Moreover, 
the multiplicity $\cM_{\cA}(t)$ of an admissible $t\in[N,2N)$ is at least 
\be\label{eq:cMAbnd}
\cM_{\cA}(t)>N^{2\gd_{\cA}-1-o(1)}
.
\ee
\end{conj}

\begin{rmk}\label{rmk:genAC}
It is now clear how to generalize \conjref{conj:McM}; the same should hold for $a_{j}$ restricted to any alphabet $\cA$, as long as $\gd_{\cA}>1/2$.
\end{rmk}

A direct attack on this conjecture
seems out of reach of current technology.
Therefore we shift our focus to the study of the arithmetic properties of the trace set, more specifically to its equidistribution along progressions, with applications to almost primes.
Our main goal is to make some 
progress in this direction. 

\subsection{Statements of the Main Theorems}\label{sec:rough}

In this subsection, we state our main theorems, though we 
defer
the precise (and somewhat technical) definitions to the next section.

For several applications, an
 important 
barometer of our understanding of
 a sequence 
 is its
 {\it level of distribution}, defined roughly as follows. 
 In our context, we wish to know 
that
the 
traces
in $\cT_{\cA}$ up to some growing parameter $N$ are equi-distributed along multiples of integers $q$, with $q$ as large as possible relative to $N$. 
That is, the quantity 
$$
\#\{t\in\cT_{\cA}: t<N,\ t\equiv0(q)\},
$$
counted with multiplicity, should be ``close'' to
$$
\frac1q\ \times\ \#\{t\in\cT_{\cA}:t<N\}
,
$$
in the sense that their difference should be much smaller than the total number of $t\in\cT_{\cA}$ up to $N$.
This proximity cannot be expected once $q$ is as large as $N$, say, but perhaps can be established with $q$ of size $N^{1/2}$ or more generally $N^{\ga}$ for some $\ga>0$. If this is the case, in an average sense, then $N^{\ga}$ is called a {\bf level of distribution} for $\cT_{\cA}$, and $\ga$ is called an {\bf exponent of distribution}.
Let us make matters a bit more precise.

Looking at traces up to $N$ counted with multiplicity is essentially the same as looking at matrices in
 the semigroup
  $\G_{\cA}$ of norm at most $N$.
Writing
$$
r_{q}(N) \ := \
\sum_{\g\in\G_{\cA}\atop\|\g\|<N}\bo_{\{\tr\g\equiv0(q)\}}
-
\frac1q
\sum_{\g\in\G_{\cA}\atop\|\g\|<N}1
$$
for the ``remainder'' terms, we will say, again roughly, that $\cT_{\cA}$ has level of distribution $\cQ$ if
\be\label{eq:levRough}
\sum_{q<\cQ}|r_{q}(N)|
=
o
\left(
\sum_{\g\in\G_{\cA}\atop\|\g\|<N}1
\right)
.
\ee
In applications it is enough to consider only square-free $q$ in the sum. 
If \eqref{eq:levRough} can be established with $\cQ$ as large as $N^{\ga}$, then we will say $\cT_{\cA}$ has exponent of distribution $\ga$.
See \secref{sec:level} for a precise definition of level and exponent of distribution.
\begin{rmk}
Note that a level and exponent of distribution is not a quantity intrinsic to $\cT_{\cA}$, but rather a function of what one can prove about $\cT_{\cA}$.
The larger this exponent, the more control one has on the distribution of $\cT_{\cA}$ on such arithmetic progressions. 
\end{rmk}
\begin{rmk}\label{rmk:AffSieve}
The set $\cT_{\cA}$ is of Affine Sieve type; see \cite{BourgainKontorovich2015a} for a definition. As such, the general Affine Sieve procedure introduced in \cite{BourgainGamburdSarnak2006, BourgainGamburdSarnak2010}, combined with
the ``expansion property'' established in \cite{BourgainGamburdSarnak2011},
 shows that $\cT_{\cA}$ has {\it some} exponent of distribution $\ga>0$, see \secref{sec:method}.
 In fact, if one replaces the known expansion by a Ramanujan-type conjecture for the spectral gap, one obtains an exponent $\ga=1/3-\vep$, see \rmkref{rmk:Raman}.
\end{rmk}

Our main goal in this paper 
is to make some partial progress towards \conjref{conj:locGlob} by 
establishing levels of distribution for $\cT_{\cA}$ beyond those available from expansion alone.
\begin{thm}\label{thm:main1}
For any small
$\eta>0$, there is an effectively computable $\gd_{0}=\gd_{0}(\eta)<1$ so that, if the dimension $\gd_{\cA}$ of the alphabet $\cA$ exceeds $\gd_{0}$, then the set $\cT_{\cA}$ 
has exponent of distribution
\be\label{eq:ga14}
\ga=\frac13-\eta.
\ee
\end{thm}

That is, 
we recover unconditionally a Ramanujan-quality exponent.
%
Applying standard sieve theory  \cite{Greaves1986}, these levels of distribution have the following immediate corollary on almost primes. 
Recall that a number is $R$-almost-prime if it has at most $R$ prime factors.

\begin{cor}\label{cor:12}
There exists an effectively computable  $\gd_{0}<1$ so that, if the dimension $\gd_{\cA}$ of the alphabet $\cA$ exceeds $\gd_{0}$, 
then the set $\cT_{\cA}$ of traces contains an infinitude of $R$-almost-primes, with $R=4$. 
\end{cor}

As an afterthought, we
explore what can be said about $R$-almost-primes, not in the set $\cT_{\cA}$ of traces, but in the set of discriminants which arise. 
To this end, recalling \eqref{eq:DMis}, we define
\be\label{eq:sDA}
\sD_{\cA}\ :=\ \{\sqf(t^{2}-4):t\in\cT_{\cA}\},
\ee 
where $\sqf(\cdot)$ denotes the square-free part. As explained in \secref{sec:Mer}, an easy consequence of Mercat's thesis \cite{Mercat2012}, combined with our work \cite{BourgainKontorovich2014} on Zaremba's Conjecture and Iwaniec's theorem \cite{Iwaniec1978}, gives
the following

\begin{thm}\label{thm:Mer}
For the alphabet $\cA=\{1,\dots,50\}$,
the set 
$\sD_{\cA}$ 
contains an infinitude of $R$-almost-primes with $R=2$. 
\end{thm}

The proof of the main \thmref{thm:main1} is based on the following result on additive energy in $\SL_2(\mathbb Z)$ of independent interest.

\begin{theorem} \label{theorem1.32}
Let $\cS_{N}=\{\gamma \in \SL_2(\mathbb Z): \Vert\gamma\Vert < N\}$.
For any sufficiently small $\kappa>0$ and $\gep>0$, there is a subset $\cS_N' \subset \cS_N$ satisfying
\be\label{eq:AddEn} 
|\cS_N \setminus \cS_N'|< N^{2-\kappa},
\ee
having additive energy
\be\label{1.34}
E(\cS_N')=\#\{(\gamma_1, \gamma_2, \gamma_3, \gamma_4)\in (\cS_N')^4: \gamma_1+ \gamma_2= \gamma_3+ \gamma_4 
\} \ll N^{4+2\kappa+\gep}.
\ee
\end{theorem}

\subsection{Organization}

In \secref{sec:state}, we give precise definitions of level and exponent of distribution, thus making unambiguous the statement of \thmref{thm:main1}.
There we also discuss the main ingredients involved in the proofs.
We spend \secref{sec:setup}  constructing the sifting sequence $\fA$, and we execute the main term analysis in \secref{sec:Main}. 
The error analysis is handled in \secref{sec:main2}, thus proving  
\thmref{thm:main1} modulo the proof of
 \thmref{theorem1.32}; the latter is postponed to \secref{sec:AddEn}. 
Finally, \thmref{thm:Mer} is proved quickly in  \secref{sec:Mer}. 

\subsection{Notation}

We use the following 
notation throughout. Set $e(x)=e^{2\pi i x}$ and $e_{q}(x)=e(\frac xq)$. We use the symbol $f\sim g$ to mean $f/g\to1$. The symbols $f\ll g$ and $f=O(g)$ 
are used 
interchangeably to mean the existence of an implied constant $C>0$ so that $f(x)\le C g(x)$ holds for all $x>C$; moreover $f\asymp g$ means $f\ll g\ll f$. 
 The letters $c$, $C$ denote  positive constants, not necessarily the same in each occurrence. 
Unless otherwise specified,  implied constants  may depend at most on $\cA$, 
which is treated as fixed.
The letter $\vep>0$ is an arbitrarily small constant, not necessarily the same at each occurrence. When it appears in an inequality, the implied constant may also depend on $\vep$ without further specification.
The symbol $\bo_{\{\cdot\}}$ is the indicator function of the event $\{\cdot\}$. The trace of a matrix $\g$ is denoted $\tr\g$.
The greatest common divisor of $n$ and $m$ is written $(n,m)$ and their least common multiple is $[n,m]$.
The function $\nu(n)$ denotes the 
number of
distinct
 prime factors of $n$. 
The cardinality of a finite set $S$ is denoted $|S|$ or $\# S$.
The transpose of a matrix $g$ is written ${}^{t}g$. When there can be no confusion, we use the shorthand $r(q)$ for $r(\mod q)$. The prime symbol $'$ in $\underset{r(q)}{\gS} {}'$ means the range of $r(\mod q)$ is restricted to $(r,q)=1$. 


\section{Levels of Distribution and 
Ingredients}\label{sec:state}

\subsection{Levels of Distribution}\label{sec:level}

In this subsection, we give  precise definitions of level and exponent of distribution.
Fix the alphabet $\cA$ and let $\cT_{\cA}$ be the set of traces of $\G_{\cA}$.
First we assume that the set of traces is {\bf primitive}, that is, 
\be\label{eq:primitive}
\gcd(\cT_{\cA})=1.
\ee
If not,%
\footnote{%
In fact, since the identity matrix has trace $2$, the set of traces is not primitive if and only if the alphabet $\cA\subset2\Z$ consists 
entirely
of even numbers (in which case the traces are all even and should be halved).%
}
then replace $\cT_{\cA}$ by $\cT_{\cA}/\gcd(\cT_{\cA})$.
Given a large parameter $N$,
let $\fA=\{a_{N}(n)\}$ be a sequence of non-negative numbers supported on $\cT_{\cA}\cap[1,N]$, and set
$$
|\fA|=\sum_{n}a_{N}(n).
$$
We require 
that
   $\fA$ is well-distributed on average over multiples of  square-free integers $\fq$. More precisely, setting 
$$
|\fA_{\fq}|:=\sum_{n\equiv0(\fq)}a_{N}(n)
,
$$
we 
insist
that
\be\label{eq:fAfq}
|\fA_{\fq}|
=
\gb(\fq)|\fA|+r(\fq)
,
\ee
where
\begin{enumerate}
\item
 the ``local density'' $\gb$ is a multiplicative function assumed to satisfy the ``linear sieve'' condition
\be\label{eq:gbEst}
\prod_{w\le p <z}
(1-\gb(p))^{-1}
\le 
C\cdot
{\log z\over \log w}
,
\ee
for some $C>1$ and any $2\le w<z$;
and
\item
the
 ``remainders''
 $r(\fq)$ are small on average, in the
 sense 
  that
 \be\label{eq:rfqLevel}
\sum_{\fq<\cQ}|r(\fq)|\ll_{K}
 {1\over (\log N)^{K}}|\fA|
,
 \ee
 for some
 $\cQ\ge1$
   and any
   $K\ge1$. That is, we ask for 
   an arbitrary power of log savings.
\end{enumerate}
If a sequence $\fA$ exists for which the conditions \eqref{eq:fAfq}--\eqref{eq:rfqLevel} hold, then we say that 
 $\cT_{\cA}$ has  a {\bf level of distribution}
$\cQ$.  If \eqref{eq:rfqLevel} can be established with $\cQ$ as large as a power, 
\be\label{eq:cQNga}
\cQ=N^{\ga},\qquad\qquad
 \ga>0,
\ee 
then we say that $\cT_{\cA}$ has an {\bf exponent of distribution} $\ga$.

\subsection{The Main Ideas}\label{sec:method}

This subsection is purely heuristic and expository.
%
%
First we
recall
 how the ``standard'' Affine Sieve procedure applies in this context, explaining Remark \ref{rmk:AffSieve}.
Since $\gd_{\cA}$ is assumed to be large, we must have $\{1,2\}\subset\cA$, whence 
for all $\fq\ge1$,
 %
 the reduction $\G_{\cA}(\mod \fq)$ is all of $\SL_{2}(\fq)$;
 cf. Remark \ref{rmk:StrongApp}.
 Initially, we could construct the sequence $\fA$ by setting
\be\label{eq:aNbad}
a_{N}(n):=\sum_{\g\in\G_{\cA}\atop\|\g\|<N}\bo_{\{\tr \g = n\}}
,
\ee
which is clearly supported on $n\in\cT_{\cA},$ $n\ll N$. Then work of Hensley \cite{Hensley1989} gives 
\be\label{eq:fAsize}
|\fA|
=
\#\{\g\in\G_{\cA}:\|\g\|<N\}
\asymp N^{2\gd_{\cA}},
\ee
and
 $|\fA_{\fq}|$ can be expressed as
\be\label{eq:fAfq1}
|\fA_{\fq}|
=
\sum_{\g\in\G_{\cA}\atop\|\g\|<N}\bo_{\{\tr \g \equiv0(\fq)\}}
=
\sum_{\g_{0}\in\SL_{2}(\fq)}\bo_{\{\tr \g_{0} \equiv0(\fq)\}}
\left[
\sum_{\g\in\G_{\cA}\atop\|\g\|<N}\bo_{\{\g \equiv \g_{0}(\fq)\}}
\right]
,
\ee
where we have decomposed the $\g$ sum into residue classes mod $\fq$. 
A theorem of Bourgain-Gamburd-Sarnak \cite{BourgainGamburdSarnak2011} in this context states very roughly (see \propref{prop:aleph} for a precise statement) that
\bea\label{eq:BGSrough}
&&
\hskip-.5in
\#\{\g\in\G_{\cA}:\|\g\|<N,\ \g\equiv \g_{0}(\fq)\}
\\
\nonumber
&&
=
\frac1{|\SL_{2}(\fq)|}
\#\{\g\in\G_{\cA}:\|\g\|<N\}
+
\text{``}O(\fq^{C}N^{2\gd-\gT})
\text{''},
\eea
for some $\gT>0$. (We reiterate that the error in \eqref{eq:BGSrough} is heuristic only;  a statement of this strength is not currently known.%
\footnote{Added in print: A power savings error now is known, and even for all $\fq$ (not just square-free) by work of Magee-Oh-Winter/Bourgain-Kontorovich-Magee \cite{MageeOhWinter2016, BourgainKontorovichMagee2015}. For
our purposes, the weaker result in \cite{BourgainGamburdSarnak2011} suffices, and in fact none of our estimates would improve (though the exposition would be slightly simpler) if we used \cite{MageeOhWinter2016, BourgainKontorovichMagee2015}
instead.} 
That said, the true statement serves the same purpose in our application.) 
This is the ``spectral gap'' or ``expander'' property of $\G_{\cA}$, and follows from a resonance-free region for the resolvent of a certain ``congruence'' transfer operator, see \secref{sec:BGS2}.

Inserting the expander property \eqref{eq:BGSrough} into $|\fA_{\fq}|$ in \eqref{eq:fAfq1} gives the desired decomposition \eqref{eq:fAfq}, with local density
$$
\gb(\fq)=\frac{1}{|\SL_{2}(\fq)|}\sum_{\g_{0}\in\SL_{2}(\fq)}\bo_{\{\tr\g_{0}\equiv0(q)\}}
,
$$
and error
\be\label{eq:rfq}
|r(\fq)|\ll \fq^{C}N^{2\gd-\gT}.
\ee
Then the local density condition \eqref{eq:gbEst} follows classically from primitivity \eqref{eq:primitive}, and, in light of \eqref{eq:fAsize}, the average error condition \eqref{eq:rfqLevel} requires (a condition weaker than)
$$
\sum_{\fq<\cQ}|r(\fq)|\ll \cQ^{C}N^{2\gd-\gT}<N^{2\gd-\vep},
$$
or 
\be\label{eq:cQbnd}
\cQ=N^{\ga}<N^{\gT/C-\vep}
.
\ee
In this way, one can prove {\it some} exponent of distribution $
\ga>0$, cf. Remark \ref{rmk:AffSieve}.
If one were to compute the numerical values of  
 the constants $C$ and $\gT$
 from the proof of \eqref{eq:BGSrough}, which would be a feat in itself,
one
would
obtain a numeric but astronomically small 
$\ga$. Our goal here is to do better.

\begin{rmk}\label{rmk:Raman}
The best one may hope to be true  is a Ramanujan-type ``square-root cancellation'' error, where $\gT=\gd$ and $C=0$ in \eqref{eq:BGSrough}, which leads to $C=2$ in \eqref{eq:rfq}
 and $C=3$ in \eqref{eq:cQbnd}. (Note that one cannot obtain an improvement here along the lines of Hong-Kontorovich \cite{HongKontorovich2015} by a better modular decomposition in \eqref{eq:fAfq1},
as it is easy to see that the trace is only stabilized mod $\fq$ by the full congruence group and not some subgroup.) The ``Ramanujan-quality'' exponent would thus be $\ga=\gT/C=\gd/3$ in \eqref{eq:cQbnd}, which approaches $1/3$ as $\gd$ approaches $1$;
this explains the claim in \rmkref{rmk:AffSieve}.
\end{rmk}

The novel technique employed here, used in some form already in \cite{BourgainKontorovich2010, BourgainKontorovich2014, BourgainKontorovich2014a,
BourgainKontorovich2015a, BourgainKontorovich2016, BourgainKontorovich2017}, is to take inspiration from Vinogradov's method, developing a ``bilinear forms'' approach, as follows.
Instead of \eqref{eq:aNbad}, let $X$ and $Y$ be two more parameters, each a power of $N$, with $XY=N$, and set (roughly)
\be\label{eq:aNbetter}
a_{N}(n)\ \text{``}:=" \
\sum_{\g\in\G\atop \|\g\|<X}
\sum_{\xi\in\G\atop\|\xi\|<Y}
\bo_{\{\tr(\g\xi)=n\}}.
\ee
This sum better encapsulates the group structure of $\G_{\cA}$, while still only being supported on the traces $\cT_{\cA}$ of $\G_{\cA}$. 
Again, this is still an oversimplification; see \secref{sec:setup} for the actual construction of $\fA$.

Instead of directly appealing to expansion as in \eqref{eq:fAfq1}, we first invoke finite abelian harmonic analysis, writing
\be\label{eq:fAfq2}
|\fA_{\fq}|=
\sum_{n\equiv0(\fq)}
a_{N}(n)
=
\sum_{n}
\left[
\frac1\fq
\sum_{r(\fq)}
e_{\fq}\left(r n\right)
\right]
a_{N}(n).
\ee
After some manipulations, we decompose our treatment according to whether $\fq$ is ``small'' or ``large''. For $\fq$ small, we apply expansion as before. For $\fq$ large, the corresponding exponential sum already has sufficient cancellation (on average over $\fq$ up to the level 
$\cQ$) that it can be treated as an error term in its entirety. It is in this range of large $\fq$ that we exploit the bilinear structure of \eqref{eq:aNbetter}. 



\subsection{Expansion}\label{sec:BGS2}


Here is a formal statement of ``expansion,'' as needed in our context.
Let $\cA\subset\N$ be 
our 
finite alphabet with 
dimension $\gd_{\cA}$ sufficiently near $1$.
As such, it must contain the sub-alphabet $\cA_{0}:=\{1,2\}\subset\cA$. This has the consequence that
for all 
$q\ge1$,
\be\label{eq:Gmodq}
\G(\mod q)\cong\SL_{2}(q)
,
\ee
cf. Remark \ref{rmk:StrongApp}. 
Furthermore, we will only require expansion for the fixed alphabet $\cA_{0}$, so as to make the expansion constants absolute, and not dependent on $\cA$; see footnote \ref{foot:rmk} on page \pageref{foot:rmk}.

To this end, let $\G_{0}\
%
\subset\SL_{2}(\Z)$
be the
  semigroup
  as  in \eqref{eq:GcAis}
  corresponding to $\cA_{0}$. 
The following proposition is proved in 
{\cite[Prop. 2.9]{BourgainKontorovich2017}}.

\begin{prop}\label{prop:aleph}
Given any $\faS\gg1$, 
there is a non-empty subset $\aleph=\aleph(\faS)\subset\G_{0}$ so that
\begin{enumerate}
\item
for all $\fa\in\aleph$, $\|\fa\|<\faS$, and
\item
for all square-free $q
$ and $\fa_{0}\in\SL_{2}(q)$, 
\be\label{eq:alephED}
\left|
{\#\{\fa\in\aleph:\fa\equiv\fa_{0}(q)\}\over |\aleph|}-{1\over |\SL_{2}(q)|}\right|
\ll
\fE(\faS;q),
\ee
where
\be\label{eq:fEbnd}
\fE(\faS;q
)
:=
\twocase{}
{e^{-c\sqrt{\log \faS}},}{if $q\le  C\log \faS,$}
{q^{C} \faS^{-\gT},}{if $q>C\log \faS.$}
\ee
\end{enumerate}
\end{prop}


\section{The Sifting Set and Initial Manipulations}\label{sec:setup}

\subsection{Construction of $\fA$}

The first goal in this subsection is to construct the appropriate sifting sequence $\fA=\{a_{N}(n)\}$. Let
$\cA\subset\N$ be our fixed alphabet with corresponding dimension $\gd_{\cA}$ near $1$, and  
let
 $\G_{\cA}$ be the semigroup in \eqref{eq:GcAis}. 
Since $\cA$ is fixed, we 
drop the subscripts, writing $\G=\G_{\cA}$ and $\gd=\gd_{\cA}$.

Let $N$ be the main growing parameter, and let
\be\label{eq:Xis}
\xiS=N^{x},\quad  \faS=N^{y}, \quad Z=N^{z},\quad x,y,z>0,\quad x+y+z=1,
\ee 
be some parameters to be chosen later;
in particular,
\be\label{eq:NXYZ}
N=\xiS \faS Z.
\ee
We think of $\xiS$ as large, $\xiS>N^{1/2}$, and $\faS$ as tiny, of size $N^{\vep}$.

Let $\aleph=\aleph(\faS)\subset\G_{0}\subset\G$ be the set constructed in Proposition \ref{prop:aleph}.
We also create certain subsets
$$
\Xi\ \subset\ \{\xi\in\Gamma: \Vert\xi\Vert<X\},\qquad
\Omega\ \subset\ \{\omega\in\Gamma: \Vert
\omega\Vert< Z\},
$$
as follows. 
Applying \thmref{theorem1.32} with $N=X$ and 
\be\label{eq:gkIs}
\kappa\ =\ 2(1-\delta)+\varepsilon
\ee (here $\vep$ is fixed, sufficiently small, and $\gd$ is sufficiently close to $1$),
gives a subset $\cS_{X}'$ of $\{\xi\in\SL_{2}(\Z): \Vert\xi\Vert<X\}$ of size
$$
|\cS_{X}\setminus\cS_{X}'|\ < \ X^{2-\gk} \ = \
X^{2\delta-\varepsilon}
,
$$
cf. \eqref{eq:AddEn}, and  having additive energy
controlled as in \eqref{1.34}.
Note that Hensley's estimate \eqref{eq:fAsize} gives
$$
\{\xi\in\Gamma: \Vert\xi\Vert<X\} \ \asymp\ X^{2\gd},
$$
which is too big to be contained in the complement $\cS_{X}\setminus\cS_{X}'$.
So defining
\be\label{4.20}
\Xi \ := \ \cS_{X}'\cap \{\xi\in\Gamma: \Vert\xi\Vert<X\}
\ee
gives a subset of $\G$ of proportional size,
\be\label {4.21}
|\Xi|\ \asymp\ X^{2\delta},
\ee
and controlled additive energy,
\be\label{4.22}  
E(\Xi) \ \le\ E(\cS_{X}') \ \ll\
X^{4+4(1-\delta)+3\varepsilon}.
\ee

We construct the set $\gW$ in the same fashion, obtaining
\be\label{eq:gWsize}
|\gW| \ \asymp \ Z^{2\gd},
\ee
and
\be\label {4.23}
E(\Omega)\ \ll\ Z^{4+4(1-\delta)+3\varepsilon}.
\ee
Note that, since $\delta$ will be taken close to 1, \eqref{4.22}, \eqref{4.23} provide a nearly optimal additive energy control.

Then we can finally define the sifting sequence $\fA=\{a_{N}(n)\}$ by:
\be\label{eq:aNdef}
a_{N}(n):=
\sum_{\xi\in\Xi}
\sum_{\fa\in\aleph}
\sum_{\gw\in\gW}
\bo_{\{n=\tr(\xi\fa\gw)\}}
.
\ee
Note that $a_{N}(n)$ is supported on $n\ll N$ by \eqref{eq:NXYZ}. 
We record from the above that
\be\label{eq:fAsizeBnd}
|\fA|=
|\Xi|\cdot
|\aleph|\cdot
|\gW|
\gg
|\aleph|
(\xiS Z)^{2\gd}
.
\ee

\subsection{Initial Manipulation}

Next for  parameters $1\ll Q_{0}<\cQ$ and any square-free $\fq<\cQ$, we decompose
\bea\nonumber
|\fA_{\fq}|
&=&
\sum_{n\equiv0(\fq)}a_{N}(n)
=
\sum_{n}
\frac1\fq
\sum_{q\mid\fq}\sideset{}{'}\sum_{r(q)}
e_{q}(rn)a_{N}(n)
\\
\label{eq:fAfqDecomp}
&=&
\cM_{\fq}+r(\fq)
,
\eea
say,
according to whether or not $q<Q_{0}$.
Here  
\be\label{eq:cMis}
\cM_{\fq}
\ :=\
\sum_{n}
\frac1\fq
\sum_{q\mid\fq\atop q<Q_{0}}\sideset{}{'}\sum_{r(q)}
e_{q}(rn)a_{N}(n)
\ee
will be treated as a ``main'' term, 
the remainder $r(\fq)$ being an error.


\section{Main Term Analysis}\label{sec:Main}

In this section, we analyze the main term, $\cM_{\fq}$, proving the following
\begin{prop}
Let $\gb$ be
the multiplicative function given at primes by
\be\label{eq:gbIs}
\gb(p)\ :=\
\frac1p
\left(
1+{\chi_{4}(p)\over p}
\right)
\left(
1-\frac1{p^{2}}
\right)^{-1}
,
\ee
where $\chi_{4}$ is the Dirichlet character mod $4$. 
There is a decomposition
\be\label{eq:cMfqD}
\cM_{\fq} \ = \ \gb(\fq)\ |\fA| + r^{(1)}(\fq)+r^{(2)}(\fq),
\ee
where 
\be\label{eq:r(1)bnd}
\sum_{\fq<\cQ}|r^{(1)}(\fq)|
\ \ll\
|\fA|
\log \cQ
\left(
{1
\over
e^{c\sqrt{\log \faS}}}
+
Q_{0}^{C}
\faS^{-\gT}
\right)
,
\ee
and
\be\label{eq:r(2)bnd}
\sum_{\fq<\cQ}|r^{(2)}(\fq)|
\ \ll\
|\fA|
{\cQ^{\vep}\over Q_{0}}
.
\ee
\end{prop}
\pf
Inserting the definition \eqref{eq:aNdef} of $a_{N}$ into \eqref{eq:cMis} gives
\beann
\cM_{\fq}
&=&
\sum_{\xi\in\Xi}
\sum_{\fa\in\aleph}
\sum_{\gw\in\gW}
\frac1\fq
\sum_{q\mid\fq\atop q<Q_{0}}\sideset{}{'}\sum_{r(q)}
e_{q}(r\tr(\xi\fa\gw))
\\
&=&
\sum_{\xi\in\Xi}
\sum_{\gw\in\gW}
\frac1\fq
\sum_{q\mid\fq\atop q<Q_{0}}\sideset{}{'}\sum_{r(q)}
\sum_{\fa_{0}\in\SL_{2}(q)}
e_{q}(r\tr(\xi\fa_{0}\gw))
\left[
\sum_{\fa\in\aleph\atop\fa\equiv\fa_{0}(q)}
1
\right]
.
\eeann

Apply \eqref{eq:alephED} to the innermost sum, giving
\beann
\cM_{\fq}
&=&
\cM_{\fq}^{(1)}
+
r^{(1)}(\fq),
\eeann
say, where
$$
\cM_{\fq}^{(1)}
\ :=\
|\fA|\
\frac1\fq
\sum_{q\mid\fq\atop q<Q_{0}}
\sideset{}{'}\sum_{r(q)}
{1\over |\SL_{2}(q)|}
\sum_{\g\in\SL_{2}(q)}
e_{q}(r\tr(\g))
,
$$
and
\beann
|r^{(1)}(\fq)|
&
\ll
&
|\fA|\
\frac1\fq
\sum_{q\mid\fq\atop q<Q_{0}}
q^{4}\
\fE(\faS;q)
.
\eeann
The error $\fE$ is as given in \eqref{eq:fEbnd}.
We estimate
\beann
\sum_{\fq<\cQ}|r^{(1)}(\fq)|
&\ll&
|\fA|
\sum_{q<Q_{0}}
q^{4}\
\fE(\faS;q)
\sum_{\fq<\cQ\atop \fq\equiv0(q)}
\frac1\fq
\\
&\ll&
|\fA|
\log\cQ
\bigg[
(\log \faS)^{
C
}
e^{-c\sqrt{\log \faS}}
+
Q_{0}^{C}
\faS^{-\gT}
\bigg]
,
\eeann
thus proving \eqref{eq:r(1)bnd}.

Returning to $\cM^{(1)}_{\fq}$, we add back in the large divisors $q\mid \fq$, writing
$$
\cM^{(1)}_{\fq}
\ =\
\cM^{(2)}_{\fq}
+
r^{(2)}(\fq)
,
$$
say, where
$$
\cM_{\fq}^{(2)}
\ :=\
|\fA|\
\frac1\fq
\sum_{q\mid\fq}
\sideset{}{'}\sum_{r(q)}
{1\over |\SL_{2}(q)|}
\sum_{\g\in\SL_{2}(q)}
e_{q}(r\tr(\g))
.
$$

Let $\rho(q)$ be the multiplicative function given at primes by
$$
\rho(p)
\ :=\
{1\over |\SL_{2}(p)|}
\sum_{\g\in\SL_{2}(p)}
\sideset{}{'}\sum_{r(p)}
e_{p}(r\tr(\g))
,
$$
so that
$$
\cM_{\fq}^{(2)}
\ =\
|\fA|\ 
\frac1\fq\prod_{p\mid\fq}\bigg(1+\rho(p)\bigg).
$$

By an elementary computation, we evaluate explicitly that
$$
\rho(p) \ = \ 
{p(p+\chi_{4}(p))\over p^{2}-1}-1,
$$
and hence 
$$
\cM_{\fq}^{(2)}
\ =
\ 
|\fA|\cdot \gb(\fq)
,
$$
with $\gb$ as given in \eqref{eq:gbIs}. 

Lastly, we deal with $r^{(2)}$. 
It is easy to see from the above that $|\rho(p)|\ll 1/p$, so $|\rho(q)|\ll q^{\vep}/q$, giving the bound
$$
|r^{(2)}(\fq)|
\ \ll\
|\fA|\
\frac1\fq
\sum_{q\mid\fq\atop q\ge Q_{0}}
\frac{q^{\vep}
}q
\ \ll\
|\fA|\
\frac{\fq^{\vep}}\fq
\frac1
{Q_{0}}
.
$$
The estimate \eqref{eq:r(2)bnd} follows immediately, completing the proof.
\epf

\begin{rmk}
Since $\faS$ in \eqref{eq:Xis} is a small power of $N$, the first error  term  in \eqref{eq:r(1)bnd} saves an arbitrary power of $\log N$, as required in \eqref{eq:rfqLevel}. For the rest of the paper, all other error terms will be power savings. In particular, setting
\be\label{eq:Q0toga0}
Q_{0}\ = \ N^{\ga_{0}}, \qquad \ga_{0}>0,
\ee
the error in \eqref{eq:r(2)bnd} is already a power savings, while the second term in \eqref{eq:r(1)bnd} requires that
\be\label{eq:ga0togT}
\ga_{0}\ < \ {y\gT\over C}.
\ee
It is here that we crucially use the expander property for $\G$, but the final level of distribution will be independent of $\gT$. 
\end{rmk}


\section{Proof of \thmref{thm:main1}} \label{sec:main2}

\subsection{Initial Manipulations}

Returning to \eqref{eq:fAfqDecomp}, it remains to control the average error term
\be\label{eq:cEdef}
\cE
\ :=\
\sum_{\fq<\cQ}|r(\fq)|
\ =\
\sum_{\fq<\cQ}
\left|
\sum_{\xi\in\Xi}
\sum_{\fa\in\aleph}
\sum_{\gw\in\gW}
\frac1\fq
\sum_{q\mid\fq\atop q\ge Q_{0}}\sideset{}{'}\sum_{r(q)}
e_{q}(r\tr(\xi\fa\gw))
\right|
,
\ee
the goal being to verify \eqref{eq:rfqLevel}.
We first massage $\cE$ into a more convenient form. 

Let $\gz(\fq):=|r(\fq)|/r(\fq)$ be the complex unit corresponding to the absolute value in \eqref{eq:cEdef}, and rearrange terms as:
$$
\cE
\ =\
\sum_{Q_{0}\le q<\cQ}
{\gz_{1}(q)\over q}
\sum_{\xi\in\Xi}
\sum_{\fa\in\aleph}
\sum_{\gw\in\gW}
\sideset{}{'}\sum_{r(q)}
e_{q}(r\tr(\xi\fa\gw))
,
$$
where we have set
$$
\gz_{1}(q)\ :=\
\sum_{\fq<\cQ/q
}
\frac{\gz(q\fq)
}\fq
.
$$
Note for future reference that
\be\label{eq:gz1}
|\gz_{1}(q)| \ll \log \cQ.
\ee

Leaving the 
special set $\aleph$ alone, we
break the $q$ sum into dyadic pieces
\be\label{eq:cEtocE1}
\cE
\ \ll\
\sum_{\fa\in\aleph}
\sum_{Q_{0}\le Q<\cQ\atop \text{dyadic}}
|\cE_{1}(Q;\fa)|
,
\ee
where we have defined
\be\label{eq:cE1Q}
\cE_{1}(Q;\fa)
\ :=\
\sum_{q\asymp Q}
{\gz_{1}(q)\over q}
\sum_{\xi\in\Xi}
\sum_{\gw\in\gW}
\sideset{}{'}\sum_{r(q)}
e_{q}(r\tr(\xi\fa\gw))
.
\ee
It remains to estimate $\cE_{1}(Q;\fa)$.
%

\subsection{Bounding $\cE_{1}$}

We claim the following estimate. 

\begin{thm}\label{prop:cEbnd2}
For any $\vep>0$, and any $1\ll 
Q_{0}<Q<\cQ<N\to\infty$, with
\be\label{eq:ZtoQ0Q1p}
Q<Z,\qquad Q^{2}=o(X)
,
\ee 
we have
\be\label{eq:cE1Qbnd2}
|\cE_{1}(Q;\fa)|
\ll
N^{1-\delta+\vep}
Q
|\gW|^{1/2}
\xiS^{2}
Z
\left[
{Q^{1/2}\over
Z^{1/2}}
+
{1\over Q^{1/8}}
\right]
.
\ee
\end{thm}
\pf
Start by applying Cauchy-Schwarz in $q,$ $r$, and the ``short'' variable $\gw$  to \eqref{eq:cE1Q}. This opens the ``long'' variable $\xi$ into a pair of such, as follows.
\beann
|\cE_{1}(Q;\fa)|^{2}
&\ll &
\left(
\sum_{q\asymp Q}
\sum_{\gw\in\gW}
\sideset{}{'}\sum_{r(q)}
{|\gz_{1}(q)|^{2}\over q^{2}}
\right)
\left(
\sum_{q\asymp Q}
\sum_{\gw\in\SL_{2}(\Z)\atop\|\gw\|<Z}
\sideset{}{'}\sum_{r(q)}
\left|
\sum_{\xi\in\Xi}
e_{q}(r\tr(\xi\fa\gw))
\right|^{2}
\right)
\\
&\ll &
|\gW|
\log^{2}\cQ
\left(
\sum_{\xi, \xi' \in \Xi} 
\left|
\sum_{q\asymp Q}
\sideset{}{'}\sum_{r(q)}
\sum_{\omega\in\Omega}
e_{q}(r\tr((\xi-\xi')\fa\gw))
\right|
\right)
.
\eeann
Collect the difference of $\xi$ and $\xi'$ into a single variable, writing  
$$
\xi-\xi'=M\in M_{2\times 2}(\Z)\cong \Z^{4}
,
$$
and setting
\be\label{eq:cNMdef}
\cN^{\Xi}_{M}(\xiS):=\sum_{\xi, \xi' \in \Xi} 
\bo_{\{M=\xi-\xi'\}}
.
\ee
In view of the additive energy bound \eqref {4.22}, we have
$$
\sum_{M\in\Z^{4}}
\cN^{\Xi}_{M}(\xiS)^{2}
\ll
\xiS^{4+\tau}
,
$$
where we have set (cf. \eqref{eq:gkIs})
$$
\tau=4(1-\delta)+3\varepsilon
.
$$
So writing
\beann
|\cE_{1}(Q;\fa)|^{2}
&\ll &
N^{\vep}
|\gW|
\sum_{M\in\Z^{4}\atop\|M\|\ll \xiS}
\cN^{\Xi}_{M}(\xiS)
\left|
\sum_{q\asymp Q}
\sideset{}{'}\sum_{r(q)}
\sum_{\omega\in \Omega} 
e_{q}(r\tr(M\fa\gw))
\right|
,
\eeann
we apply Cauchy-Schwarz in the $M$ variable, giving
\beann
|\cE_{1}(Q;\fa)|^{4}
&\ll &
N^{\vep}
|\gW|^{2}
X^{4+\tau}
\sum_{M\in\Z^{4}}
\Psi\left({M\over \xiS}\right)
\left|
\sum_{q\asymp Q}
\sideset{}{'}\sum_{r(q)}
\sum_{\omega\in\Omega} 
e_{q}(r\tr(M\fa\gw))
\right|^{2}
\\
&\ll &
N^{\vep}
X^{\tau}
|\gW|^{2}
X^{4}
\sum_{q,q'\asymp Q}
\sideset{}{'}\sum_{r(q)}
\sideset{}{'}\sum_{r'(q')}
\sum_{\omega, \omega' \in\Omega} 
\\
&&
\hskip1in
\sum_{M\in\Z^{4}}
\Psi\left({M\over \xiS}\right)
e\left(M\cdot\left(\frac rq\fa\gw-\frac{r'}{q'}\fa\gw'\right)\right)
\\
&\ll &
X^{\tau}
Q^{4}
|\gW|^{2}
\xiS^{8}
\sum_{q,q'\asymp Q}
\sideset{}{'}\sum_{r(q)}
\sideset{}{'}\sum_{r'(q')}
\sum_{\omega, \omega'\in \Omega}  
\bo_{\{\|\frac{q' r\fa\gw-{r'}{q}\fa\gw'}{qq'}\|< \frac1\xiS\}}
.
\eeann
Here we have inserted a suitable bump function $\Psi$ and applied Poisson summation to the sum on $M\in\Z^{4}$. 
Assuming $qq'\ll Q^{2}=o(X)$ as in \eqref{eq:ZtoQ0Q1p}, the innermost condition implies
$$
q' r\gw\equiv q{r'}\gw'
(\mod qq')
.
$$
Taking determinants gives
$$
(q' r)^{2}\equiv({r'}{q})^{2}
(\mod qq')
,
$$
whence reducing mod $q$ gives
$$
(q' r)^{2}\equiv0
(\mod q)
.
$$
But this implies $(q')^{2}\equiv0(\mod q)$, since $(r,q)=1$.  Because $q$ is square-free, we have thus forced $q'\equiv0(\mod q)$. By symmetry, we also have $q\equiv 0(\mod q')$, and hence
$$
q=q',\qquad
r\equiv ur'(\mod q),\qquad\text{and}\qquad
\gw\equiv u\gw'(\mod q),
$$
where $u^{2}\equiv1(q)$; 
note that
there are at most $2^{\nu(q)}\ll N^{\vep}$ such $u$'s. We then have
\beann
|\cE_{1}(Q;\fa)|^{4}
&\ll &
N^{\vep}
X^{\tau}
|\gW|^{2}
\xiS^{8}
\sum_{q\asymp Q}
\sum_{u^{2}\equiv1(q)}
\sideset{}{'}\sum_{r,r'(q)\atop r\equiv ur'(q)}
\sum_{\omega, \omega' \in\Omega}  
\bo_{\{\gw\equiv u\gw'(\mod q)\}}
.
\eeann
We dispose of $u$ in the last summation via Cauchy-Schwarz:
\beann
\sum_{\omega, \omega'\in \Omega} 
\bo_{\{\gw\equiv u\gw'(\mod q)\}}
&=&
\sum_{\g\in\SL_{2}(q)}
\left[
\sum_{\omega\in\Omega} 
\bo_{\{\gw\equiv\g(q)\}}
\right]
\left[
\sum_{\omega'\in\Omega} 
\bo_{\{u\gw'\equiv\g(q)\}}
\right]
\\
&\le&
\left(
\sum_{\g\in\SL_{2}(q)}
\left[
\sum_{\omega\in\Omega} 
\bo_{\{\gw\equiv\g(q)\}}
\right]^{2}\
\right)^{1/2}
\\
&&\times
\left(
\sum_{\g\in\SL_{2}(q)}
\left[
\sum_{\omega'\in\Omega} 
\bo_{\{u\gw'\equiv\g(q)\}}
\right]^{2}\
\right)^{1/2}
\\
&=&
\sum_{\omega, \omega'\in \Omega} 
\bo_{\{\gw\equiv \gw'(\mod q)\}}
,
\eeann
since $(u,q)=1$. 
Applying this estimate gives
\beann
|\cE_{1}(Q;\fa)|^{4}
&\ll &
N^{\vep}
X^{\tau}
Q\ 
|\gW|^{2}
\xiS^{8}
\sum_{q\asymp Q}
\sum_{\omega, \omega'\in\Omega} 
\bo_{\{\gw\equiv \gw'(\mod q)\}}
\\
&= &
N^{\vep}
X^{\tau}
Q\ 
|\gW|^{2}
\xiS^{8}
\sum_{q\asymp Q}
\sum_{M\in\Z^{4}\atop M\equiv0(q)}
\cN^{\gW}_{M}(Z)
,
\eeann
where we have now set
$$
\cN^{\gW}_{M}(Z)
:=
\sum_{\omega, \omega'\in\Omega} 
\bo_{\{\gw- \gw'=M\}}
.
$$

We first isolate the $M=0$ term, writing
\be\label{eq:cE14}
|\cE_{1}(Q;\fa)|^{4}
\ll
N^{\vep}
X^{\tau}
Q \ 
|\gW|^{2}
\xiS^{8}
\left(
Q
Z^{2}
+
\cE_{2}
\right)
,
\ee
say, where
$$
\cE_{2}:=
\sum_{q\asymp Q}
\sum_{M\in\Z^{4}\atop M\equiv0(q), M\neq0}
\cN^{\gW}_{M}(Z)
.
$$

Applying the Cauchy-Schwarz inequality in the variables $q$ and $M$ and recalling \eqref{eq:ZtoQ0Q1p}, i.e. $Q<Z$ gives
\be
\begin{aligned}
\mathcal E^2_2& \ll Q\Big(\frac ZQ\Big)^4 \sum_{q\asymp Q} \sum_{M\in\mathbb Z^4, M\not= 0\atop M\equiv 0(q)} \mathcal N^{\gW}_M(Z)^2\\
&\ll \frac {Z^4}{Q^3} \sum_{M\in\mathbb Z^4, M\not= 0} \mathcal N^{\gW}_M(Z)^2\sum_{q|M} 1 \ll \frac {Z^4}{Q^3} N^\varepsilon 
\sum_M\mathcal N^{\gW}_M(Z)^2\\
\label{eq:cE2Bnd}
&\ll \frac {Z^{8+\tau}}{Q^3} N^\varepsilon
\end{aligned}
\ee
where this time the additive energy bound \eqref {4.23}  was used.

Combining \eqref{eq:cE2Bnd} with \eqref{eq:cE14} gives \eqref{eq:cE1Qbnd2},
as claimed.
\epf

\subsection{Proof of \thmref{thm:main1}}

We give a sketch, as the details are very similar to \cite[\S6]{BourgainKontorovich2017}.
Inserting \eqref {eq:cE1Qbnd2} 
in \eqref{eq:cEtocE1} and recalling \eqref{eq:gWsize} and \eqref{eq:fAsizeBnd}  gives
\be\label{5.7}
\cE\ll N^{3(1-\delta)+\varepsilon} |\fA|
\left({\cQ^{1/2}\over Z^{1/2}} + {1\over Q_{0}^{1/8}}\right).
\ee
With $\gd$ very near $1$, the second term in  \eqref{5.7} is a  power saving as long as $Q_{0}$ is some tiny power.\footnote{\label{foot:rmk}It is here that we crucially need expansion to a {\it fixed} alphabet, so we can move $\gd$ near $1$ without affecting $Q_{0}$; cf. \secref{sec:BGS2}.}
Recalling \eqref{eq:Q0toga0}, \eqref{eq:ga0togT}, this  requires $Y=N^{y}$ with $y>0$ an arbitrary small fixed exponent (taking $\delta$ accordingly close to 1).
Writing  $\cQ=N^{\ga}$, $X=N^{x}$, $Z=N^{z}$ with 
$$
1=x+y+z\approx x+z
$$
the first term in \eqref{5.7} gives a power savings provided  $z>\ga+6(1-\delta+\varepsilon)$, while \eqref{eq:ZtoQ0Q1p} is satisfied for $x>2\ga$. 
Hence, taking $\delta =\delta (\eta)$ close enough to 1, we reach exponent of distribution $\alpha>\frac 13 -\eta$.
This proves \thmref{thm:main1}.


\section{Proof of \thmref{theorem1.32}}\label{sec:AddEn}

Recall our notation 
$$
\cS_{N}=\{\g\in\SL_{2}(\Z):\|\g\|<N\},
$$ 
and for an integer matrix $M\in M_{2\times2}(\Z)= \Z^{4} $, set
$$
\cN_{N}(M) \ : = \
\sum_{\g,\g'\in\cS_{N}}\bo_{\{\g+\g'=M\}}
.
$$
This is the number of representations of $M$ as a sum of two elements of $\cS_{N}$.
The additive energy of $\cS_{N}$ is the square of the $L^{2}$ norm of $\cN_{N}$:
$$
E(\cS_{N}) \ := \
\#\{\g_{1},\g_{2},\g_{3},\g_{4}\in\cS_{N} \ : \
\g_{1}+\g_{2}=\g_{3}+\g_{4}
\}
\ = \
\sum_{M\in \Z^{4}}\cN_{N}(M)^{2}.
$$
\begin{rmk}
\begin{itemize}
\item
When $M=\mattwos 0000$, we clearly have 
$$
\cN_{N}(M)=|\cS_{N}|\asymp N^{2},
$$ 
so this one term already contributes $N^{4}$ to the energy. 
\item 
If $M=\mattwos1000$, then the number of representations, $\cN_{N}(M)$, is of order $N$, since for all $n\ll N$,
\be\label{eq:Mmany}
\mattwo n1{-1}0-\mattwo {n-1}1{-1}0\ = \ M.
\ee
Such $M$'s have determinant zero, and we show below that this is the key feature allowing \eqref{eq:Mmany}. There are 
about $N^{2}$ such $M=\mattwos abcd$ (since once we fix $|a|,|d|\le N$,  there are only $N^{\vep}$ values for $b,c$ satisfying $ad=bc$), and each contributes $\cN_{N}(M)^{2}\asymp N^{2}$ to the energy, for a net contribution of around $N^{4}$.
\item
Generically, we expect $\cN_{N}(M)$ to be of size $N^{\vep}$, which again would contribute $N^{4}$ to the energy, since there are this many generic $M$.
\end{itemize}
\end{rmk}

It is thus reasonable to make the following
\begin{conj}\label{conj:addEn}
The additive energy of $\SL_{2}(\Z)$ is as small as possible,
$$
E(\cS_{N})\ll N^{4+\vep}.
$$
\end{conj}

While we are not able to prove this full conjecture, \thmref{theorem1.32} will be a sufficiently strong substitute in our applications.

Our first goal is to understand straight lines in $G:=\SL_{2}(\R)$; these are responsible for the behavior in \eqref{eq:Mmany}.
\begin{lem}
Two matrices $A, B\in G$ lie on a line in $G$; that is,  for all $t\in\R$, 
\be\label{eq:ABline}
tA+(1-t)B\ \in \ G,
\ee
if and only if
$$
\det(A-B)=0.
$$
\end{lem}
\pf
We use the elementary formula:
\be\label{eq:detForm}
\det(X+Y)  \ = \
\det(X)
+\det(Y)
+\det(X)
\tr(X^{-1}\cdot Y).
\ee
By \eqref{eq:ABline}, we have
$$
1 \ = \ \det(tA+(1-t)B)  
\ = \
t^{2} 
+(1-t)^{2} 
+t^{2} 
\tr(t^{-1}A^{-1}\cdot (1-t)B),
$$
which simplifies to
$$
(\tr(A^{-1}B) -2)t(t-1)
 \ =  \ 
 0.
$$
This equation holds for all $t$ if and only if 
$\tr(A^{-1}B)=2$.
Now using
%
 \eqref{eq:detForm} again gives
$$
\det(A-B)  \ = \
2-
\tr(A^{-1}\cdot B)
,
$$
from which the claim follows.
\epf

\begin{rmk}
This explains the phenomenon observed in \eqref{eq:Mmany}. Indeed, let $\g_{t}=B+t(A-B)=tA+(1-t)B$ as in \eqref{eq:ABline}. So if $A,B$ and $t$ are all integral, and hence $\g_{t}\in\SL_{2}(\Z)$, then $M=A-B$ has many  representations as $\g_{t}-\g_{t-1}=M$. 
\end{rmk}

Next recall the Cartan decomposition $G=KAK$, that is, any $g\in G$ can be expressed as 
$$
g=k_{\gt}a_{t}k_{\vf},
$$ 
where
$$
k_{\gt}\ = \ 
\mattwo{\cos\gt}{\sin\gt}{-\sin\gt}{\cos\gt}
,\qquad
\text{ and }\qquad
a_{t} \ = \
\mattwo{e^{t}}{}{}{e^{-t}},\ t\ge0.
$$

\begin{lem}\label{lem:slab}
There is a small absolute constant $c>0$ with the following property.
Let $A_{1},A_{2}\in\cS_{N}
$, 
and write 
$$
A_{j} \ = \ 
k_{\gt_{j}}a_{t_{j}}k_{\vf_{j}}
.
$$
Assume that
$$
|\gt_{1}-\gt_{2}|<\frac cN, \quad
|t_{1}-t_{2}|<c, \quad\text{ and } \quad
|\vf_{1}-\vf_{2}|<\frac cN
.
$$
Then
$$
\det(A_{1}-A_{2})=0.
$$
\end{lem}
\pf
Using \eqref{eq:detForm} yet again gives
$$
\det(A_{1}-A_{2})
\ = \ 
2-
\tr(
A_{1}
^{-1}\cdot 
A_{2}
)
\ = \
2-
\tr(
k_{\gt_{2}-\gt_{1}}a_{t_{2}}k_{\vf_{2}-\vf_{1}}
a_{-t_{1}}
)
.
$$
Write $\gt'=\gt_{2}-\gt_{1}$ and $\vf'=\vf_{2}-\vf_{1}$. One then computes:
\beann
\tr(
k_{\gt'}a_{t_{2}}k_{\vf'}
a_{-t_{1}}
)
& =&
\cos(\gt')\cos(\vf')
(e^{t_{1}-t_{2}}+e^{t_{2}-t_{1}})
\\
&&\hskip1in
-
\sin(\gt')\sin(\vf')
(e^{t_{1}+t_{2}}+e^{-t_{1}-t_{2}})
\\
&=&
(1+O(c^{2}N^{-2}))
(2+O(c))
\\
&&\hskip1in
-
O(c^{2}N^{-2})
O(
N^{2})
\\
&=&
2+O(c)
,
\eeann
since $c<1$. Thus we have $\det(A_{1}-A_{2})=O(c)$, but since this is an integer, we can take $c$ small enough (independent of $N$) to force the desired conclusion.
\epf

In light of this lemma, we 
restrict our attention to
``sub-slabs'' of $\cS_{N}$ 
in which the $\gt$ and $\vf$ parameters are restricted to intervals of length $c/N$, and $t$ to an interval of length $c\asymp1$. That is, for each such triplet $\gT_{\ga},T_{\ga},\Phi_{\ga}$ of intervals of length $
|\gT_{\ga}|=|\Phi_{\ga}|=c/N$,
 $|T_{\ga}|=c$, 
 set
$$
\cS_{N,\ga} \ := \ \{A=k_{\gt}a_{t}k_{\vf}\in \cS_{N} :(\gt,t,\vf)\in\gT_{\ga}\times T_{\ga}\times \Phi_{\ga}\}.
$$
We will require $O(N^{2})$ such to cover $\cS_{N}$:
$$
\cS_{N} = \bigsqcup_{\ga}\cS_{N,\ga}.
$$
The previous lemma tells us that the elements of a slab $\cS_{N,\ga}$, if any, snap to a single affine line in $\SL_{2}(\R)\subset\R^{4}$. If we can make the ``slope'' of this line large, then the points will be far-spaced, implying that there can be only very few of them in a single slab. Let us make this precise.
\\

Assume that $A, B\in\cS_{N,\ga}$; then $\det(A-B)=0$ by \lemref{lem:slab}. Set $M=B-A$; this is the ``slope.'' Since $\det M=0$, we can write $ M$ as:
\be\label {6.19}
 M =\begin{pmatrix} rv_1& sv_1\\ rv_2& sv_2\end{pmatrix}
,
\ee
 in which $v=  \begin{pmatrix} v_1 \\ v_2 \end{pmatrix}$ 
 is a primitive vector in  $ \mathbb Z^2$, and $r, s\in\mathbb Z$.
 Observe moreover that
 $$
\det B \ = \ \det(A+ M) \ = \
\det A+\det  M +\det A\cdot\tr(A^{-1}\cdot M),
 $$
 where we again used \eqref{eq:detForm}. Hence
\be\label{eq:trAxi}
 \tr(A^{-1}\cdot M)=0,
\ee
which, on writing $A=\mattwos abcd$, gives 
\be
\label{6.20}
drv_1- brv_2 -csv_1+ asv_2 =0.
\ee
Multiplying \eqref{6.20} by $a$ and using $ad=1-bc$ gives
$$
rv_{1}
=
bcrv_1+abrv_2 +acsv_1- a^{2}sv_2
=
(cv_1+av_2)(br-as)
,
$$
which implies that
\be\label{6.21}
|cv_1+av_2| \ |br-as|
\leq |rv_1| \leq\| M\|.
\ee
This equation gives our requisite lower bound on the slope. We can finally define the desired subset $\cS_{N}'$ of $\cS_{N}$.

\begin{Def}
Given small parameters $\gk>0$ and $\gep>0$, set 
$$
\gk_{1}=\gk+\gep,
$$
and define
\be\label{6.22}
\begin{aligned}
\cS_N'\ = \  &\Bigg\{ A=\begin{pmatrix} a&b\\ c&d\end{pmatrix} \in SL_2(\mathbb Z): \Vert A\Vert\leq N \text{ and}\\
&  \ \min_{\substack{k\in\mathbb Z^2\\ 0<|k|<N^{1-\kappa_1}}}
|k|\cdot\min\big(|k_1 a+k_2b|, |k_1 a+k_2c|\big) > N^{1-\kappa_1}\Bigg\}.
\end{aligned}
\ee
\end{Def}

\begin{lem}
We have
\be\label{6.23}
|\cS_N\setminus \cS_N'|\ \ll\ N^{2-\kappa_1+\varepsilon}
.
\ee
\end{lem}
Thus \eqref {eq:AddEn}  holds, since $\kappa <\kappa_1$. 

\pf
We will bound the number of $A\in\SL_{2}(\Z)$ with $\|A\|<N$ and for which there exists some $0<|k|<N^{1-\gk_{1}}$ so that
$$
 |k_1 a+k_2b| \ \le\ N^{1-\gk_{1}}/|k|
 .
$$
Break $|k|$ dyadically into $|k|\asymp K$, with $K<N^{1-\gk_{1}}$. We may assume $|a|\ge |b|>0$ (since $b=0$ gives only parabolic elements), so $(a,b)=1$, and break $|a|$ dyadically into $|a|\asymp M$, with $M\le N$.
Write 
$$
y\ = \ k_{1}a+k_{2}b,
$$ 
so that 
$$
|y| \ \ll \ {N^{1-\gk_{1}}\over K}.
$$
Thus we have
\be\label{eq:cSNbnd1}
|\cS_{N}\setminus\cS_{N}'|
\ \ll \
\sum_{K<N^{1-\gk_{1}}\atop\text{dyadic}}
\sum_{M\le N\atop\text{dyadic}}
\sum_{k\in\Z^{2}\atop |k|\asymp K}
\sum_{y\in\Z\atop |y|\ll {N^{1-\gk_{1}}\over K}}
\sum_{(a,b)=1, |b|\le |a|\asymp M\atop y=k_{1}a+k_{2}b}
\sum_{(c,d)\in\Z^{2}\atop |(c,d)|<N}
\bo_{\{ad-bc=1\}}
.
\ee
We will estimate this expression from the inside out.

Once $(a,b)$ is determined, the equation $ad-bc=1$ has a unique solution $(c_{0},d_{0})$ with $0\le d_{0}<|a|$ and $0\le c_{0}<|b|$, and all other solutions $(c,d)$ satisfy
$
c\equiv c_{0}(\mod a).
$
Of course, once $c$ is determined, so is $d$, and hence the last summation of \eqref{eq:cSNbnd1} contributes
$$
\ll \ 1 + {N\over M}.
$$

To determine $(a,b)$, we  handle separately the cases $y=0$ or not. 
If the former, that is, $k_{1}a=-k_{2}b$, we first choose $k_{1}$ and $a$ (for which there are $KM$ choices), and then there are $N^{\vep}$ choices for $k_{2}$ and $b$. 
This gives a net contribution to \eqref{eq:cSNbnd1} of
$$
\ll\ 
\sum_{K<N^{1-\gk_{1}}\atop\text{dyadic}}
\sum_{M\le N\atop\text{dyadic}}
KM
N^{\vep}
\left(
1+\frac NM
\right)
\ \ll \
N^{2-\gk_{1}+\vep}
.
$$

Now  we assume $y\neq0$. We want to exploit $k_{2}b= y-k_{1}a$, so we must handle separately whether $y-k_{1}a=0$ or not. 
In the case $y=k_{1}a$, we choose $y$ first (with at most $N^{1-\gk_{1}}/K$ possible values), whence there are $N^{\vep}$ choices for $k_{1}$ and $a$ (recall $y\neq0$).
Since  $b\neq0$, the condition $k_{2}b= y-k_{1}a=0$ implies that $k_{2}=0$; then we are free to choose $b$ (in $M$ choices). In total, the contribution of this case to \eqref{eq:cSNbnd1} is
$$
\ll \
\sum_{K<N^{1-\gk_{1}}\atop\text{dyadic}}
\sum_{M\le N\atop\text{dyadic}}
 {N^{1-\gk_{1}}\over K}
N^{\vep}
M
\left(
1+\frac NM
\right)
\ \ll \
N^{2-\gk_{1}+\vep}
.
$$

Lastly we consider the case $y\neq0$ and $y\neq k_{1}a$. We fix $y$, $k_{1}$, and $a$, with $N^{1-\gk-1}/K \cdot K\cdot M$ choices. Then there are $N^{\vep}$ choices for $k_{2}$ and $b$ satisfying $k_{2}b=y-k_{1}a$, giving a final contribution of
$$
 \ll \
\sum_{K<N^{1-\gk_{1}}\atop\text{dyadic}}
\sum_{M\le N\atop\text{dyadic}}
{N^{1-\gk_{1}}\over K}
KMN^{\vep}
\left(
1+\frac NM
\right)
\ \ll \
N^{2-\gk_{1}+\vep}
$$
to \eqref{eq:cSNbnd1}. This completes the proof.
\epf

Now suppose $A,B\in\cS_{N,\ga}\cap \cS_{N}'$ and $ M=B-A\neq0$. We claim that $\|M\|\ge N^{1-\gk_{1}}$. Assume by contradiction that $\|M\|<N^{1-\gk_{1}}$. Then $ M$ is of the form \eqref{6.19} satisfying also \eqref{eq:trAxi}, and moreover the vectors $(v_{2},v_{1}), (-s,r)$ are of size at most $N^{1-\gk_{1}}$, so enter in the definition of $\cS_{N}'$ in \eqref{6.22}.
By this definition, we have that 
$$
\| M\|\ \ge \ |cv_1+av_2| \ |br-as|
\ > \
{N^{1-\gk_{1}}\over |v|}{N^{1-\gk_{1}}\over |(r,s)|}
\ \gg \
{N^{2-2\gk_{1}}\over \| M\|}
,
$$
so
$$
\|A-B\|  \ = \ \| M\| \ \gg \ N^{1-\gk_{1}}
 \ = \ N^{1-\gk-\gep}
.
$$
That is, any such $A$ and $B$, while all  lying on a single line, are also very much spaced apart. This immediately gives the following
\begin{cor}
$$
|\cS_{N,\ga}\cap \cS_{N}'| \ \ll \ N^{\gk+\gep},
$$
and hence the additive energy of each such slab is
\be\label {6.25}
E(\cS_{N,\ga}\cap \cS_{N}') \ \ll \ N^{2\gk+\gep}|\cS_{N,\ga}\cap \cS_{N}'|.
\ee
\end{cor}

To put the various slabs and their energies together, we appeal to the recent resolution \cite{BourgainDemeter2015} of the $L^{2}$-decoupling conjecture; more precisely, we require the version for ``generalized cones''  proved in \cite[Theorem 1.1]{Oh2016}, in the following form.

\begin{thm}\label{thm:L2}
The energy of $\cS_{N}'$ is controlled by that of its slabs by:
\be\label{6.9}
E(\cS_N') \ \ll \  N^{2+\varepsilon} \left\{\sum_\alpha E(\cS_N'\cap \cS_{N,\ga})\right\}.
\ee
\end{thm}
Before explaining how \eqref{6.9} follows from \cite{Oh2016}, let us first observe that we now have proved \eqref{1.34} and hence \thmref{theorem1.32}. 
\pf[Proof of \thmref{theorem1.32} assuming \thmref{thm:L2}]
Indeed, combining \eqref{6.9} with \eqref{6.25} gives
$$
E(\cS_N') 
\ \ll \  N^{2+\varepsilon} 
\sum_\alpha N^{2\gk}|\cS_{N,\ga}\cap \cS_{N}'|
\ = \  N^{2+2\gk+\varepsilon+\gep} 
|\cS_{N}'|
\ \ll \  N^{4+2\gk+\gep} .
$$
\epf
\begin{rmk}
While we believe \conjref{conj:addEn} for the full set $\cS_{N}$, our use of $L^{2}$-decoupling makes it absolutely essential to excise certain regions of $\cS_{N}$, leaving only $\cS_{N}'$. Indeed, the region responsible for \eqref{eq:Mmany}, that is,  having $\gt\approx0, \ \vf\approx 0$, already contributes  $N^{3}$ to the energy of its slabs;
together with an extra loss of size $N^{2}$ on the right side of \eqref{6.9}, this would  give an unacceptable net contribution of $N^{5}$ to the energy.
\end{rmk}

It remains to establish \thmref{thm:L2}.
\pf
The linear map 
$$
\mattwos abcd\mapsto \left(\frac {a+d}2, \frac {a-d}2, \frac {b-c}2, \frac {b+c}2\right)
$$
identifies $SL_2(\mathbb R)$ and the hyperboloid 
$$
\tilde C:x^2-y^2+z^2-w^2=1
.
$$
Thus 
$$
\frac 1N(\tilde C \cap \{\xi\in\mathbb R^4: \Vert\xi\Vert\le N\})
$$ 
is the set 
$$
\{\xi=(x, y, z, w): \Vert \xi\Vert\le 1
\text{
and }x^2-y^2+z^2 -w^2 =\frac 1{N^2}\}
.
$$
This lies within a $\frac 1{N^2}$-neighborhood
of 
$$
\{\xi\in C: \Vert\xi\Vert\le 1\},
$$
where $C$ is the  cone 
$$
C: x^2-y^2 +z^2 -w^2=0 . 
$$

We will use the decoupling theorem obtained in \cite[Therem 1.1]{Oh2016} to the conical surface
taking $p=4$, $\delta \asymp \frac 1{N^2}$.
Denote $C_\delta$ a $\delta$-neighborhood of $C$.
In the parametrization
\be\label{6.3}
x= r\cos \theta, y= r \cos \psi, z= r \sin \theta, w= r\sin\psi
\ee
the slab-decomposition of $C_\delta$ is obtained by restricting $(\gt, \psi)\in I_\alpha \times J_\alpha$, $I_\alpha, J_\alpha$ intervals of size $\delta^{\frac 12}\asymp \frac 1N$.
Denoting $\{\tau_\alpha\}$ the corresponding slabs, we obtain therefore (denoting by $\hat f$ the Fourier transform of $f$ on $\mathbb R^4$).

\begin{lemma}\label{6.4}
Let $\text {\,supp\,} \hat f \subset C_\delta \cap [\Vert \xi\Vert\asymp 1]$ and $f_\alpha =(\hat f|_{\tau_\alpha})^\vee$ the Fourier restriction
of $f$ to $\tau_\alpha$. Then
\be\label {6.5}
\Vert f\Vert_{L^4(B_{N^2})} \ \ll \ N^{\frac 12+\varepsilon} \left(\sum_\alpha \Vert f_\alpha\Vert^4_{L^4(B_{N^2})}\right)^{1/4}.
\ee
\end{lemma}

Here we have denoted
$$
\Vert f\Vert_{L^p(B_K)} \ = \ \left(\int_{\mathbb R^4} |f(x)|^p \Big(1+\frac {|x|}K\Big)^{-100} dx\right)^{\frac 1p}
$$
(i.e. the $L^p$-norm of $f$ on an appropriately `smoothed' ball \break $\{x\in\mathbb R^4:|x|<K\}$).

According to the previous discussion, 
$$ 
\tilde C\cap[\Vert\xi\Vert\asymp N]  \ \subset \ C_{\frac 1N}\cap [\Vert\xi\Vert \asymp N].
$$
Rescaling \eqref{6.5} implies therefore

\begin{lemma}\label {6.6}
Let $\text{supp\,} \hat f \subset\tilde C_{\frac 1N} \cap [\Vert\xi\Vert\asymp N]$ and $f_\alpha =(\hat f|_{N\tau_\alpha})^\vee$.
Then
\be\label{6.7}
\Vert f\Vert_{L^4(B_N)} \ \ll \ N^{\frac 12+\varepsilon}\left(\sum_{\alpha} \Vert f_\alpha\Vert^4_{L^4(B_N)}\right)^{\frac 14}.
\ee
\end{lemma}
\medskip 

To prove \thmref{thm:L2}, we apply \lemref{6.6}  in the discretized setting, taking
\be\label{6.8}
f(x) =\sum_{\xi\in \cS_N'} e(x.\xi).
\ee


Note that \eqref{6.8} is a 1-periodic function in $x$ and hence
$$
N^{-4}\Vert f\Vert^4_{L^4(B_N)}\ \asymp\ \int_{[0, 1]^4}|f(x)|^4 dx\ =  \ E(\cS_N').
$$
Thus \eqref{6.7} allows us to bound the additive energy of $\cS_N'$ in terms of 
the additive energies of its subsets $\cS_N' \cap N\tau_\alpha$; it remains to reinterpret intervals of the form $N\tau_{\ga}$ as $S_{N,\ga}$.

As discussed previously, we identify  $G=\SL_2(\mathbb R)$ and $\tilde C$, so that
 $\cS_N'$ may be thought of as a subset of $\tilde C$.
Writing $A\in G$ as 
\be\label{6.11}
A=k_{u}a_{t}k_{v}
\ee
we have
$$
A=e^{t}  \begin{pmatrix} \cos u \cos v&  \cos u \sin v \\ -\sin u  \cos v&-\sin u  \sin v\end{pmatrix} + O\left(\frac 1N\right)
$$
and hence, by \eqref{6.3}
$$
\begin{cases}
r \cos \theta = \frac 12 e^{t} \cos (u+v) +O(\frac 1N)\\
r \cos\psi = \frac 12 e^{t} \cos (u-v) +O(\frac 1N)\\
r\sin \theta = \frac 12  e^{t} \sin (u+v)+ O(\frac 1N)\\
r \sin \psi = \frac 12  e^{t}\sin (v-u)+ O(\frac 1N).
\end{cases}
$$
Thus
$$
\begin{cases}
 e^{t}=2r+O(\frac 1N)\\
\cos\theta= \cos (u+v)+O(\frac 1{N^2})\\
\cos\psi = \cos(u-v)+O(\frac 1{N^2})\\
\sin\theta =\sin (u+v)+O(\frac 1{N^2})\\
\sin\psi=\sin(v-u)+O(\frac 1{N^2}).
\end{cases}
$$
Restriction of $(\theta, \psi)\in I_\alpha \times J_\alpha$ corresponds therefore to a restriction of $(u, v)\in I_\alpha'\times J_\alpha '$, with $I_\alpha', J_\alpha'$ size
$\frac 1N$ intervals.
Hence $\cS_N' \cap N\tau_\alpha$ corresponds to $\cS_N' \cap \cS_{N,\ga}$, as needed.
\epf


\section{Proof of \thmref{thm:Mer}}\label{sec:Mer}

Recall from \eqref{eq:sDA} that
$\sD_{\cA}$ is the set of discriminants which arise from the alphabet $\cA$. Set $\cA=\{1,\dots,A\}$ with $A=50$.
 In his thesis, Mercat connects the Arithmetic Chaos Conjecture with Zaremba's, by proving the following 
\begin{thm}[\cite{Mercat2012}]
If the reduced rational $m/n$ has all partial quotients bounded by $A$, and if the denominator $n$ arises as a solution to the Pellian equation $n^{2}-\gD r^{2}=\pm1$, then $\Q[\sqrt \gD]\cap\fC_{\cA}$ is non-empty. 
\end{thm}
In fact, he exhibits a periodic continued fraction in $\Q[\sqrt \gD]$ via an explicit construction involving the partial quotients of $m/n$.

With his theorem, we can now sketch a  

\pf[Proof of \thmref{thm:Mer}]
Iwaniec's theorem \cite{Iwaniec1978} states that the number of $n$ up to $N$ with $\gD=n^{2}+1$ having at most $2$ prime factors is at least $CN/\log N$. Taking the alphabet $\cA=\{1,\dots,50\}$ in \cite{BourgainKontorovich2014},  the exceptional set is of order much smaller than $N/\log N$, and hence $100\%$ of such denominators $n$ have a coprime numerator $m$ with $m/n$ having all partial quotients bounded by $A=50$. Clearly setting $r=1$ gives a solution to $n^{2}-\gD r^{2}=-1$, whence $\gD\in\sD_{\cA}$ by Mercat's theorem.
\epf


\section*{Acknowledgments} 

It is our pleasure to thank
Curt McMullen
for many detailed comments and suggestions on an earlier version of this paper, and
  Tim Browning, Zeev Rudnick, and Peter Sarnak for illuminating conversations. 
  Thanks also to Michael Rubinstein for numerics related to
\conjref{conj:addEn}.  
The second-named author would like to thank the hospitality of the IAS, where much of this work was carried out.

\bibliographystyle{amsplain}

\begin{thebibliography}{MOW16}

\bibitem[Art24]{Artin1924}
Emil Artin.
\newblock Ein mechanisches system mit quasiergodischen bahnen.
\newblock {\em Abh. Math. Sem. Univ. Hamburg}, 3(1):170--175, 1924.

\bibitem[BD15]{BourgainDemeter2015}
Jean Bourgain and Ciprian Demeter.
\newblock The proof of the {$l^2$} decoupling conjecture.
\newblock {\em Ann. of Math. (2)}, 182(1):351--389, 2015.

\bibitem[BGS06]{BourgainGamburdSarnak2006}
Jean Bourgain, Alex Gamburd, and Peter Sarnak.
\newblock Sieving and expanders.
\newblock {\em C. R. Math. Acad. Sci. Paris}, 343(3):155--159, 2006.

\bibitem[BGS10]{BourgainGamburdSarnak2010}
Jean Bourgain, Alex Gamburd, and Peter Sarnak.
\newblock Affine linear sieve, expanders, and sum-product.
\newblock {\em Invent. Math.}, 179(3):559--644, 2010.

\bibitem[BGS11]{BourgainGamburdSarnak2011}
J.~Bourgain, A.~Gamburd, and P.~Sarnak.
\newblock Generalization of {S}elberg's 3/16th theorem and affine sieve.
\newblock {\em Acta Math}, 207:255--290, 2011.

\bibitem[BK10]{BourgainKontorovich2010}
J.~Bourgain and A.~Kontorovich.
\newblock On representations of integers in thin subgroups of {SL}$(2,{{\bf
  {Z}}})$.
\newblock {\em GAFA}, 20(5):1144--1174, 2010.

\bibitem[BK14a]{BourgainKontorovich2014}
J.~Bourgain and A.~Kontorovich.
\newblock On {Z}aremba's conjecture.
\newblock {\em Annals Math.}, 180(1):137--196, 2014.

\bibitem[BK14b]{BourgainKontorovich2014a}
Jean Bourgain and Alex Kontorovich.
\newblock On the local-global conjecture for integral {A}pollonian gaskets.
\newblock {\em Invent. Math.}, 196(3):589--650, 2014.

\bibitem[BK15]{BourgainKontorovich2015a}
Jean Bourgain and Alex Kontorovich.
\newblock The {A}ffine {S}ieve {B}eyond {E}xpansion {I}: {T}hin {H}ypotenuses.
\newblock {\em Int. Math. Res. Not. IMRN}, (19):9175--9205, 2015.

\bibitem[BK16]{BourgainKontorovich2016}
J.~Bourgain and A.~Kontorovich.
\newblock Beyond {E}xpansion {I}{I}{I}: {R}eciprocal {G}eodesics, 2016.
\newblock  {\tt arXiv:1610.07260}.

\bibitem[BK17]{BourgainKontorovich2017}
J.~Bourgain and A.~Kontorovich.
\newblock Beyond expansion {I}{I}: Low-lying fundamental geodesics.
\newblock {\em J. Eur. Math. Soc. (JEMS)}, 19(5):1331--1359, 2017.



\bibitem[BKM15]{BourgainKontorovichMagee2015}
J.~Bourgain, A.~Kontorovich, and M.~Magee.
\newblock Thermodynamic expansion to arbitrary moduli, 2015.
\newblock To appear, {\em Crelle's Journal} {\tt arXiv:1507.07993}.

\bibitem[FK14]{FrolenkovKan2014}
Dmitrii~A. Frolenkov and Igor~D. Kan.
\newblock A strengthening of a theorem of {B}ourgain-{K}ontorovich {II}.
\newblock {\em Mosc. J. Comb. Number Theory}, 4(1):78--117, 2014.

\bibitem[Goo41]{Good1941}
I.~J. Good.
\newblock The fractional dimensional theory of continued fractions.
\newblock {\em Proc. Cambridge Philos. Soc.}, 37:199--228, 1941.

\bibitem[Gre86]{Greaves1986}
G.~Greaves.
\newblock The weighted linear sieve and {S}elberg's {$\lambda^2$}-method.
\newblock {\em Acta Arith.}, 47(1):71--96, 1986.

\bibitem[Hen89]{Hensley1989}
Doug Hensley.
\newblock The distribution of badly approximable numbers and continuants with
  bounded digits.
\newblock In {\em Th\'eorie des nombres ({Q}uebec, {PQ}, 1987)}, pages
  371--385. de Gruyter, Berlin, 1989.

\bibitem[Hen92]{Hensley1992}
Doug Hensley.
\newblock Continued fraction {C}antor sets, {H}ausdorff dimension, and
  functional analysis.
\newblock {\em J. Number Theory}, 40(3):336--358, 1992.

\bibitem[Hen96]{Hensley1996}
Douglas Hensley.
\newblock A polynomial time algorithm for the {H}ausdorff dimension of
  continued fraction {C}antor sets.
\newblock {\em J. Number Theory}, 58(1):9--45, 1996.

\bibitem[HK15]{HongKontorovich2015}
Jiuzu Hong and Alex Kontorovich.
\newblock Almost prime coordinates for anisotropic and thin pythagorean orbits.
\newblock {\em Israel J. Math.}, 209(1):397--420, 2015.

\bibitem[Hua15]{Huang2015}
ShinnYih Huang.
\newblock An improvement to {Z}aremba's conjecture.
\newblock {\em Geom. Funct. Anal.}, 25(3):860--914, 2015.

\bibitem[Hum16]{Humbert1916}
G.~Humbert.
\newblock Sur les fractions continues ordinaires et les formes quadratiques
  binaires ind{\'e}finies.
\newblock {\em Journal de math{\'e}matiques pures et appliqu{\'e}es 7e
  s{\'e}rie}, 2:104--154, 1916.

\bibitem[Iwa78]{Iwaniec1978}
Henryk Iwaniec.
\newblock Almost-primes represented by quadratic polynomials.
\newblock {\em Invent. Math.}, 47:171--188, 1978.

\bibitem[Jen04]{Jenkinson2004}
Oliver Jenkinson.
\newblock On the density of {H}ausdorff dimensions of bounded type continued
  fraction sets: the {T}exan conjecture.
\newblock {\em Stoch. Dyn.}, 4(1):63--76, 2004.

\bibitem[McM09]{McMullen2009}
Curtis~T. McMullen.
\newblock Uniformly {D}iophantine numbers in a fixed real quadratic field.
\newblock {\em Compos. Math.}, 145(4):827--844, 2009.

\bibitem[McM12]{McMullenNotes}
C.~McMullen.
\newblock Dynamics of units and packing constants of ideals, 2012.
\newblock Online lecture notes,
  \url{http://www.math.harvard.edu/~ctm/expositions/home/text/papers/cf/slides/slides.pdf}.

\bibitem[Mer12]{Mercat2012}
P.~Mercat.
\newblock Construction de fractions continues p\'eriodiques uniform\'ement
  born\'ees, 2012.
\newblock To appear, {\it J. Th\'eor. Nombres Bordeaux}.

\bibitem[MOW16]{MageeOhWinter2016}
M.~Magee, H.~Oh, and D.~Winter.
\newblock Uniform congruence counting for {S}chottky semigroups in
  {S}{L}$(2,{Z})$, 2016.
\newblock {\tt arXiv:1601.03705}.

\bibitem[MVW84]{MatthewsVasersteinWeisfeiler1984}
C.~Matthews, L.~Vaserstein, and B.~Weisfeiler.
\newblock Congruence properties of {Z}ariski-dense subgroups.
\newblock {\em Proc. London Math. Soc}, 48:514--532, 1984.

\bibitem[Oh16]{Oh2016}
Changkeun Oh.
\newblock Remarks on {W}olff's inequality for hypersurfaces, 2016.
\newblock {\tt arXiv:1602.05861v2}.

\bibitem[Ser85]{Series1985}
Caroline Series.
\newblock The modular surface and continued fractions.
\newblock {\em J. London Math. Soc. (2)}, 31(1):69--80, 1985.

\bibitem[Wil80]{Wilson1980}
S.~M.~J. Wilson.
\newblock Limit points in the {L}agrange spectrum of a quadratic field.
\newblock {\em Bull. Soc. Math. France}, 108:137--141, 1980.

\bibitem[Woo78]{Woods1978}
A.~C. Woods.
\newblock The {M}arkoff spectrum of an algebraic number field.
\newblock {\em J. Austral. Math. Soc. Ser. A}, 25(4):486--488, 1978.

\end{thebibliography}


\begin{dajauthors}
\begin{authorinfo}[pgom]
Jean Bourgain\\
School of Mathematics\\
Institute for Advanced Study\\
Princeton, NJ\\
bourgain\imageat{}ias\imagedot{}edu \\
  \url{https://www.math.ias.edu/people/faculty/bourgain}
\end{authorinfo}
\begin{authorinfo}[johan]
Alex Kontorovich\\
Department of Mathematics\\ 
Rutgers University\\
New Brunswick, NJ\\
alex.kontorovich\imageat{}rutgers\imagedot{}edu \\
  \url{http://sites.math.rutgers.edu/~alexk/}
\end{authorinfo}
\end{dajauthors}

\end{document}